\documentclass[10pt]{amsart}
\usepackage{amssymb,latexsym,amscd}

\numberwithin{equation}{section}

\newtheorem{Thm}{Theorem}[section]
\newtheorem{Cor}[Thm]{Corollary}
\newtheorem{Prop}[Thm]{Proposition}
\newtheorem{Lem}[Thm]{Lemma}
\theoremstyle{definition}
\newtheorem{Def}[Thm]{Definition}
\newtheorem{Ex}[Thm]{Example}
\newtheorem{Rmk}[Thm]{Remark}

\newtheorem{Conj}[Thm]{Conjecture}
\date{}

\newcommand{\LL}{{\mathcal{L}}} 
  
\newcommand{\II}{{\mathcal{I}}}   
\newcommand{\OO}{{\mathcal{O}}}
\newcommand{\PP}{{\mathbb{P}}}
\newcommand{\LP}{{\mathcal{L}_{n,d}}(-\sum_{i=1}^hm_iP_i)}
\newcommand{\ih}{i=1, \dots, h}
\newcommand{\bPP}{\tilde{\PP}}

\begin{document}

\title{Special effect varieties in higher dimension}
\author{Cristiano Bocci}

\address{Cristiano Bocci\\ 
Dipartimento di Matematica, Universit\`a di Milano\\
Via Cesare Saldini 50, 20133 Milano, Italy}
\email{Cristiano.Bocci@unimi.it}
\thanks{This reseach was partially supported by GNSAGA of INdAM (Italy).}
\begin{abstract} 
Here we introduce the concept of special effect varieties in higher
dimension and we generalize to $\PP^n$, $n\geq 3$, the two conjectures given in \cite{sev-1}
for the planar case. Finally, we propose some examples on the product of projective spaces
and we show how these results fit with the ones of Catalisano,
Geramita and Gimigliano.
\end{abstract}
\dedicatory{In memory of my grandfather Annibale}
\maketitle


\section{Introduction}

Let $\LL_{n,d}:=|\OO_{\PP^n}(d)|$ be the  complete linear system of divisors of degree
$d$ in $\PP^n$. 
Fix points $P_1,\dots, P_h$ on $\PP^n$ in general position and positive 
integers $m_1, \dots ,m_h$. We denote by $\LP$ the subsystem of $\LL$ given 
by all divisors having multiplicity at least $m_i$ at $P_i$, $i=1, \dots, h$. 
Since a point of multiplicity $m$ imposes 
$\left(\begin{smallmatrix} m+n-1 \\ n \end{smallmatrix}\right)$ conditions
we can define the {\bf virtual dimension} $\nu$ of the system $\LP$ as
$$\nu(\LL_{n,d}(-\sum_{i=1}^hm_iP_i)):=\binom{d+n}{n} -1 -\sum_{i=1}^h \binom{m_i+n-1}{n}.$$
The virtual dimension can be computed on the blow-up $\pi:\bPP^n \to \PP^n$ at the points $P_1, \dots, P_h$.
In fact, let $E_i$, $\ih$ be the exceptional divisors corresponding 
to the blow-up of the points $P_i, \ih$ and denote by $H$ the
pull-back of a general hyperplane of $\PP^n$ via $\pi$, in such a way
we can write the strict transform of the system 
$\LL:=\LL_{n,d}(\sum_{i=1}^hm_iP_i)$ 
as $\tilde{\LL}=|dH-\sum_{i=1}^hm_iE_i|$. 
It is an easy application of the (generalized) Riemann-Roch theorem to observe that
\begin{equation}\label{chi}
\nu(\LL)=\chi(\tilde{\LL})-1=h^0(\bPP^n,\tilde{\LL})-h^1(\bPP^n,\tilde{\LL})-1.
\end{equation}
We then define the 
{\bf expected dimension} $\epsilon$ of $\LP$ as
$$\epsilon(\LP):=\mbox{max}\{\nu(\LP),-1\}.$$
Since the conditions imposed by the multiple points $m_iP_i$ could be 
dependent, in general we have
$$
\dim(\LP)\geq\epsilon(\LP)
$$
We say that a system $\LP$ is {\bf special} if strict inequality holds,
 otherwise $\LP$ is said to be {\bf non-special}.

Starting with the case $X=\PP^2$, we have some precise 
conjectures about the characterization of special linear systems and 
a rich series of results on the conjectures.
The main Conjectures are the following.
  
\begin{Conj}[(SC) B. Segre, 1961]\label{SC} If a linear system 
of plane curves with general multiple base points 
$\LL_{2,d}(-\sum_{i=1}^hm_iP_i)$ is special, then its general 
member is non-reduced, i.e. the linear system has, 
according to Bertini's theorem, some multiple fixed component.
\end{Conj}

\begin{Conj}[(HHC) Harbourne-Hirschowitz, 1989]\label{HHC} A linear system 
of plane curves with general multiple base points 
$\LL:=\LL_{2,d}(-\sum_{i=1}^hm_iP_i)$ is special if and only if is $(-1)-$special, 
i.e. its strict transform on the blow-up along the points $P_1, \dots,
P_h$ splits as 
$\tilde{\LL}=\sum_{i=1}^kN_iC_i+\tilde{{\mathcal{M}}}$
where the $C_i$, $i=1, \dots ,k$, are $(-1)-$curves such that 
$C_i\cdot \tilde{\LL}=-N_i<0$,  $\nu(\tilde{{\mathcal{M}}})\geq 0$ and 
there is at least one index $j$ such that $N_j>2$.
\end{Conj}

In \cite{CiMi3} C. Ciliberto and R. Miranda proved 
that the Harbourne--Hirschowitz and Segre Conjectures are equivalent.

In \cite{sev-1} the concepts of $\alpha-$special effect curve and
$h^1-$special effect curve are introduced and two new conjectures are
proposed (see Definitions \ref{numspec} and \ref{cohomspec} for ``numerically'' and ``cohomologically'' special).

\begin{Conj}[(NSEC) ``Numerical Special Effect'' Conjecture]\label{NSEC}
A linear system of plane curves $\LL_{2,d}(-\sum_{i=1}^hm_iP_i)$ 
with general multiple base points is special if and only if it is numerically special.
\end{Conj}

\begin{Conj}[(CSEC) ``Cohomological Special Effect'' Conjecture]\label{CSEC}
A linear system of plane curves  $\LL:=\LL_{2,d}(-\sum_{i=1}^hm_iP_i)$ 
with general multiple base points is special if and only if it is cohomologically special.
\end{Conj}

The main result in \cite{sev-1} is the following
\begin{Thm}\label{equiv}
Conjectures {\rm (SC)}, {\rm (HHC)}, {\rm (NSEC)} and {\rm (CSEC)} are equivalent.
\end{Thm}

When we pass to $\PP^n$, $n\geq 3$, very little is known 
about special linear systems.
One of the most important result
is the classification of the homogeneous special systems for double points:
\label{listaspeciali}
\begin{Thm}[Alexander--Hirschowitz, 1996]\label{P^nspecial}
The system ${\mathcal{L}}_{n,d}(2^h)$ is non-special unless:
$$\begin{array}{cccccc}
n  & any & 2 & 3 & 4 & 4 \cr
d  & 2   & 4 & 4 & 4 & 3 \cr 
h  & 2, \dots, n & 5 & 9 & 14 & 7 \cr
\end{array}$$
\end{Thm} 

Continuing with $\PP^n$, $n\geq3$
we can notice that there is not a precise conjecture.
Although the Segre Conjecture can be generalized
in every ambient variety using the statement concerning 
$H^1\not=0$ (see, for example, \cite{Bocci}, \cite{BoMi} or \cite{CiMi3})
there is nothing that characterizes the special systems
from a geometric point of view as, for example,
in the case of $(-1)-$curves in $\PP^2$.

A worthy goal would be ``find a conjecture (C) in $\PP^n$,
[or in a generic variety $X$] such that, when we read (C) in $\PP^2$,
(C) is equivalent to the Segre (\ref{SC})
and Harbourne--Hirschowitz (\ref{HHC}) Conjectures''.

A first conjecture in this direction was given in \cite{Ciliberto}
where the speciality of a system in $\PP^n$ was related to the
existence of rational curves in the base locus with particular properties on their normal bundle.
Recently, Laface and Ugaglia found a counterexample to this conjecture
(see \cite{LafUga}). They showed that the linear system
$\LL:=\LL_{3,9}(-6P_0-\sum_{i=1}^84P_1)$ 
in $\PP^3$ is special and 
the only curve
contained in its base locus has genus 2. 

As already observed, Theorem \ref{equiv} assure us that both
{\it Numerical Special Effect Conjecture}
and {\it Cohomological Special Effect Conjecture}
are potential candidates for the above-mentioned goal.

In Sections \ref{alphahigh} and \ref{h1high} we generalize the
special effect curves in $\PP^2$ to special effect varieties in $\PP^n$.
The main goal of these sections is to prove that
the Conjectures hold for every special system listed in Theorem
\ref{P^nspecial}.

In Section \ref{sevexs} we present some interesting examples of
special effect varieties. In particular we show that the special system
in the Laface--Ugaglia example is both numerically and
cohomologically special.

In Section \ref{alphacross} we give some interesting
evidence 
about a possible generalization of the Numerical 
Conjecture to linear systems in the product of projective spaces. 
In particular we observe how our results fit with similar results
given in \cite{CaGeGi2},
\cite{CaGeGi3},  \cite{CaGeGi4}.

Whenever not otherwise specified, we work over the field ${\mathbb{C}}$. 

\medskip
\noindent{\bf Acknowledgements.} The author would like to thank
L. Chiantini, C. Ciliberto, A.V. Geramita, A. Gimigliano, A. Laface and
R. Miranda for useful discussions.


\section{$\alpha-$Special effect varieties in $\PP^n$, $n\geq 3$}\label{alphahigh}

Let $\LL:=\LL_{n,d}(-\sum_{i=1}^hm_iP_i)$ 
be an effective linear system on $\PP^n$. When we blow up $\PP^n$ at
the points $P_i$, $i=1,\dots, h$, we can write
\begin{equation}\label{nunu}
\nu(\LL):=\chi(\tilde{\LL})-1=h^0(\bPP^n,\tilde{\LL})-h^1(\bPP^n,\tilde{\LL})-1.
\end{equation}

Let $Y\subset \PP^n$ be a variety with codim$(Y,\PP^n)\geq 1$ passing
through some of the points $P_1, \dots, P_h$. We define $\LL-Y:=\tilde{\LL}\otimes \II_{\tilde{Y}}$
The main question we could pose is
if we can use the $\chi$ of a certain invertible sheaf
as in the case of multiple points to compute $\nu(\LL-\alpha Y)$.
For example let $\bPP^n$ be the blow-up of $\PP^n$ at the points $P_1, \dots, P_h$
and let $\LL':=\tilde{\LL}$  be the strict transform of $\LL$.
After that, we blow up $\tilde{\PP}^n$ along $\tilde{Y}$
and compute $\chi(\tilde{\LL'}-\alpha R)$,
where $R$ is the exceptional divisor $\PP({\mathcal{N}}_{\tilde{Y}|\tilde{\PP}^n})$.
We can ask if $\nu(\LL-\alpha Y)=\chi(\tilde{\LL'}-\alpha R)-1$.
Unfortunately this method does not work for every $Y$.
This is due to the fact that after the two blow-ups
some extra-generators can appear in $H^i(\tilde{\LL}'-\alpha R)$ for $i\geq 2$.
Then it can happen that $h^0(\tilde{\LL}'-\alpha R)=0$, but $\chi(\tilde{\LL}'-\alpha R)>0$, 
that is the system is empty although we expect it to be nonempty.
Thus we define the virtual dimension of a system $\LL-Y$ as
\[
\nu(\LL-Y)=h^0(\LL\otimes \II_Y)-h^1(\LL\otimes \II_Y)-1.
\]
By (\ref{nunu}) we see that this definition fits with the standard
one. Moreover, it fits with the results of Laface and Ugaglia in \cite{P3}.

We observe that, in this way, the speciality of the system is given by the non-vanishing of $H^1(\LL\otimes \II_Y)$,
that is exactly what we expected by the generalization of the Segre
Conjecture (see \cite{Ciliberto}).

\begin{Def}\label{aSEPhigh} 
Let $\LL$ and $P_1, \dots P_h$ as above. 
An irreducible variety $Y$ has the 
$\alpha-${\bf special effect property} for $\LL$ on $\PP^n$ 
if there exist positive integer $\alpha, c_{j_1}, \dots c_{j_s}$, such that
\begin{itemize}
\item[(i)]  $Y$ contains the point $P_{j_i}$ with multiplicity at
  least $c_{j_i}$  for $j=1, \dots, s$, where $P_{j_i}\in \{ P_1, \dots, P_h \}$; 
\item[(ii)] $\nu(\LL-\alpha Y)> \nu(\LL)$.
\end{itemize}
and, if  codim$(Y,\PP^n)=1$, we require
$\alpha e \leq d$ and
$1 \leq \alpha \leq \mbox{min}\{\lceil \frac{m_{j_i}}{c_{j_i}}\rceil,
i=1, \dots, s\}$, where $e:=\deg(Y)$.
Moreover we require that $\alpha$ is the maximum admissible value for
the $\alpha-$special effect property and, if $\beta > \alpha$ then
$\nu(\LL-\beta Y)<\nu(\LL-\alpha Y)$.
\end{Def}

\begin{Def}\label{aSEVhigh} 
Let $\LL$ and $P_1, \dots P_h$ be as above. An irreducible variety $Y$ is 
an $\alpha-${\bf special effect variety} for $\LL$ on $\PP^n$ if $Y$ has the 
$\alpha-$special effect property for $\LL$ and moreover $\nu(\LL-\alpha Y)\geq 0$.
\end{Def}

\begin{Def}\label{aSEChigh}
Let $\LL$ be a system as above.
Fix a sequence of (not necessarily distinct)
irreducible varieties $Y_1, \dots Y_t$,
Suppose further that
\begin{itemize}
\item[(1)] $Y_j$ has the $\alpha_j-$special effect property 
for $\LL-\sum_{i=1}^{j-1}\alpha_iY_i$, for $j=1,\dots, t$,
\item[(2)] $\nu(\LL-\sum_{i=1}^t\alpha_i Y_i)\geq 0$.
\end{itemize}
Then we call both $X:=\sum_{i=1}^\alpha Y_i$ and $\{ Y_1, \dots, Y_r\}$
an $(\alpha_1, \dots, \alpha_r)-${\bf special effect configuration for
$\LL$}.
\end{Def}

\begin{Rmk} \rm It is possible to use the $\chi$ of a certain invertible
sheaf in several situations, for example
in the case of homogeneous systems, i.e. when
$m=m_1= \dots=m_h$ with $Y$ smooth, irreducible, $c_1=\dots=c_h=1$
and with $\alpha=m$, i.e. $\alpha$ exhausts the multiplicity at the points.
In this situation we blow up $\PP^n$ along $Y$
obtaining an exceptional divisor $R$;
then condition $(ii)$ becomes
\begin{equation}\label{specialineq}
\chi(dH-\alpha R)>\chi(dH-\sum_{i=1}^hmE_i)
\end{equation}
where the $\chi$ on the left side is taken on $X=Bl_{Y}(\PP^n)$
while the $\chi$ on the right side is taken on $X'=Bl_{\{P_i\}}(\PP^n)$,
i.e. the blow-up of $\PP^n$ at $P_1, \dots P_h$.
\end{Rmk}

\begin{Rmk}\rm
In general we refer to conditions $(ii)$ of Definition \ref{aSEPhigh},
condition $(2)$ of Definition \ref{aSEChigh}
and formula (\ref{specialineq}) as the {\bf special inequality}.  
\end{Rmk}

Let $X$ be an $\alpha-$special effect variety or an $(\alpha_1, \dots,
\alpha_r)-$special effect configuration for a system $\LL$.
Then $X$ forces $\LL$ to be special. 
In fact, one has  
$$\dim(\LL)\geq\dim(\LL-X )\geq\nu(\LL-X)>\nu(\LL)$$
and, together with condition $\nu(\LL-X)\geq0$, one has $\dim(\LL)>\epsilon(\LL)$.

These facts permit us to define a particular kind of speciality. 

\begin{Def}\label{numspec}
A special system arising from the existence of an 
$\alpha-$special effect variety (or an $(\alpha_1, \dots, \alpha_r)-$special 
effect configuration) is called {\bf Numerically Special}.
\end{Def}

Finally, we can state the same conjecture as in the planar case:

\begin{Conj}[(NSEC) ``Numerical Special Effect'' Conjecture]\label{HNSEC}
A linear system $\LL_{n,d}(-\sum_{i=1}^hm_iP_i)$ 
with general multiple base points is special if and only if it is numerically special.
\end{Conj}

We restrict now our attention to $2-$special effect varieties in $\PP^n$, $n\geq3$
for the homogeneous case $\LL_{n,d}(2^h)$.
In particular we consider as a special effect variety
respectively a smooth divisor $Y=dH$,
a linear space $Y=\PP^s$, $1\leq s \leq n$
and a rational normal curve $C_n \subset \PP^n$,
i.e. the image of $\PP^1$ under the $n-$Veronese embedding.

\subsection{Hypersurfaces in $\PP^n$}
If $Y$ is a smooth hypersurface of degree $e$ passing through $P_1,
\dots, P_h$, then
the conditions for $Y$ to be a $2-$special effect variety for $\LL:=\LL_{n,d}(2^h)$ become
\begin{align}
&\binom{e+n}{ n} -1 \geq h\label{con1div}\\
&\binom{d-2e+n}{n} >  \binom{d+n}{n}- h(n+1).\label{con3div} \\
&\binom{d-2e+n}{n} \geq 1, \mbox{ i.e. } d \geq 2e \label{con2div}
\end{align}
We have the following

\begin{Prop}\label{dH}
Let $Y$ be a smooth hypersurface passing through $P_1, \dots, P_h$,
$h\geq n$. Then $Y$ is a $2-$special effect variety for
$\LL_{n,d}(2^h)$, $n \geq 2$ when
\begin{itemize}
\item[i)] $Y=\PP^{n-1}$, for $\LL_{n,2}(2^n)$ $\forall n\geq 2$;
\item[ii)] $Y=Conic \subset \PP^2$,  for $\LL_{2,4}(2^5)$;
\item[iii)] $Y=Quadric \subset \PP^3$,  for $\LL_{3,4}(2^9)$;
\item[iv)] $Y=Quadric \subset \PP^4$, for  $\LL_{4,4}(2^{14})$. 
\end{itemize}
\end{Prop}

\begin{proof}
From conditions (\ref{con1div}) and (\ref{con3div})
we obtain the following bounds on $h$:
$$\binom{e+n}{n} - 1 \geq h > \frac{1}{n+1}
\left[ \binom{d+n}{n}-\binom{d-2e+n}{n} \right].$$ 
\noindent
Then our special effect variety exists if $d\geq 2e$ and if
\begin{equation}\label{phiPr}
\begin{split}
\varphi(d,e,n):&= \binom{d+n}{n}-\binom{d-2e+n}{n}  
- (n+1) \binom{e+n}{n} +n+1
\end{split}
\end{equation}
is negative. By Pascal's triangle for binomials we have
\[
\varphi(d,e,n)=\sum_{i=1}^{2e}\binom{d+n-i}{n-1} 
- (n+1) \binom{e+n}{n} + (n+1).
\]
Thus the function $\varphi(d,e,n)$ is increasing monotone in $d$
and we can fix our attention in the case $\varphi(2e,e,n)$ and,
eventually increase the value of $d$ up to the first $d_0$
such that $\varphi(d_0,e,n)\geq 0$. Then our equation becomes
\begin{equation}\label{phid=2e}
\varphi(2e,e,n):=\binom{2e+n}{n}+ n 
- (n+1) \binom{e+n}{n}.
\end{equation}
\noindent
{\bf Step 1a:} $e=1, d=2$
\begin{equation*}
\begin{split}
\varphi(2,1,n)&=\binom{2+n}{n}+ n - (n+1) \binom{1+n}{n} = \\
&=\frac{1}{2}(n+1)(n+2)+n-(n+1)^2=\frac{1}{2}n(1-n)<0 \quad \forall n\geq 2.
\end{split}
\end{equation*}
\noindent
{\bf Step 1b:} $e=1, d\geq 3$
\begin{equation*}
\begin{split}
\varphi(d,1,n)&=\binom{d+n-1}{n-1}+ \binom{d+n-2}{n-1}
+n+1- (n+1) \binom{n+1}{n}  \geq \\
&\geq \frac{1}{6}n(n+1)(n+2)+\frac{1}{2}n(n+1)+n+1-(n+1)^2=\\
&=\frac{n}{6}(n^2-1 )>0 \quad \forall n\geq 2.
\end{split}
\end{equation*}
From conditions (\ref{con1div}) and $h\geq n$ (to avoid degenerate
cases) one has  $h=n$. Then the
first case of the Proposition follows.
We will obtain again this result in Proposition \ref{syslin}.

\noindent
{\bf Step 2a:} $e=2, d=4$
\begin{equation*}
\begin{split}
\varphi(4,2,n)&=\binom{4+n}{n}+ n 
- (n+1) \binom{2+n}{n} = \\
&=\frac{n}{24}(n^3-2n^2-13n+14)
\end{split}
\end{equation*}
and we found
$$\begin{array}{cl}
\varphi(4,2,2)<0, & h=5; \\
\varphi(4,2,3)<0, & h=9; \\
\varphi(4,2,4)<0, & h=14; \\
\varphi(4,2,n)>0, & \mbox{ for }  n\geq 5.
\end{array}$$
\noindent
{\bf Step 2b:} $e=2, d\geq 5$
\begin{equation*}
\begin{split}
\varphi(d,2,n)= & \sum_{i=1}^4 \binom{d+n-i}{n-1}+n+1 
- (n+1) \binom{2+n}{n} \geq \\
&\geq \frac{1}{120}n(n^4+15n^3+25n^2-15n-26)>0 \quad \forall n\geq 2.
\end{split}
\end{equation*}
Then, also cases $(ii)$, $(iii)$ and $(iv)$ of the Proposition are proved.

\noindent
{\bf Step 3:} at this point, we can reduce the study of $\varphi(d,e,n)$
with the condition $d \geq 2e \geq 6$.
But, in this case, $\varphi(d,e,n)\geq 0$ $\forall n\geq 3$,
as stated in the next lemma, and this concludes the proof.
\end{proof}

\begin{Lem}[Numerical Lemma 1]\label{LN1}
Let $\varphi(d,e,n)$ be defined as in (\ref{phiPr}). If  $d \geq 2e \geq 6$, 
then $\varphi(d,e,n)\geq 0$, $\forall n\geq 3$. 
\end{Lem}

\begin{proof}
Since $\varphi(d,e,n)$ is non-decreasing in $d$, we fix our attention 
on the minimal value $d=2e$. Thus we write
\begin{equation*}
\begin{split}
 \varphi(2e,e,n) & =\frac{(n+1)\dots(n+2e)}{2e!}+n-\frac{(n+1)(n+1)\dots(n+e)}{e!}= \\
 &=n + \frac{(n+1)\dots(n+e)}{e!}\lbrack A(e) \rbrack
\end{split}
\end{equation*}
where  
$$A(e):=\frac{(n+e+1)\dots(n+2e)}{(e+1)\dots(2e)}-n-1.$$
\noindent 
We show that $A\geq 0$, $\forall e \geq 3$ and $\forall n\geq3$.

\noindent
{\bf Claim:} $A(e+1) > A(e)$

\noindent
We have
$$A(e+1):=\frac{(n+e+2)\dots(n+2e+2)}{(e+2)\dots(2e+2)}-n-1.$$
$$A(e):=\frac{(n+e+1)\dots(n+2e)}{(e+1)\dots(2e)}-n-1.$$
\noindent
Thus we compute
$$A(e+1)-A(e) = \frac{(n+e+2)\dots(n+2e)}{(e+1)\dots(2e+2)}\lbrack n^2(e+1)+n(e+1)\rbrack > 0 $$
\noindent
and the claim follows.\\
\noindent
If we consider $e=3$ we obtain
$$A(3)=\frac{1}{120}\lbrack(n+4)(n+5)(n+6)\rbrack-n-1
=\frac{1}{120}n(n^2+15n-46)>0 \quad \forall n\geq 3$$
and so, for the claim, we have
$$\varphi(d,e,n)\geq 0 \quad \forall d\geq 6, \forall e \geq 3, \forall n \geq 3.$$ 
\end{proof}

\begin{Rmk}\label{alpha1dH}\rm
It is easy to see that, under the hypothesis on $Y$ as in Proposition \ref{dH}, the case $\alpha=1$
does not give any new special effect hypersurfaces 
other than the ones in Proposition \ref{dH}.
In fact the conditions for $Y$ to be a $1-$special effect variety
for a system $\LL_{n,d}(2^h)$ are
\begin{align*}
&\binom{e+n}{n} -1 \geq h\\
&\binom{d-e+n}{n} -h >  \binom{d+n}{n}- h(n+1),\\
&\binom{d-e+n}{n}-h \geq 1.
\end{align*}
Then our special effect variety can exist if
\begin{equation}\label{phiPra1}
\begin{split}
\psi(d,2,n):&= \binom{d+n}{n}-\binom{d-e+n}{n}  
- n\binom{e+n}{n} +n = \\
&= \sum_{i=1}^{e}\binom{d+n-i}{n-1} 
- n \binom{e+n}{n} + n <0. 
\end{split}
\end{equation}
Once again, we can consider the minimal value $d=2e$:
\begin{equation}
\psi(2e,e,n):=\binom{2e+n}{n}+ n 
- (n+1) \binom{e+n}{n}<0
\end{equation}
Since $\psi(2e,e,n)$ is equal to $\varphi(2e,e,n)$ in (\ref{phid=2e})
the case $\alpha=1$ does not produce any new examples of special effect hypersurfaces.
\end{Rmk}

\begin{Rmk}\rm
The argument of the proof of Proposition \ref{dH}
can be used succesfully when the system $\LL$ is homogeneous.
In general we use the equations given by the numerical
speciality 
to construct a function $\varphi$
such that our problem of the existence of an $\alpha-$special effect variety
can become a pure combinatorial problem.
The function $\varphi$ can change depending on the data of the variety $Y$,
the system $\LL$ and the ambient variety $X$. 
\end{Rmk}

\subsection{Linear Spaces in $\PP^n$}
Let $Y$ be a linear space $\PP^s \subset \PP^n$ with $1\leq s \leq n-1$;
by changing the coordinates,
we can suppose that $Y$ is defined by $x_0=\dots=x_{n-s-1}=0$.
It is not difficult to verify that the expected dimension of $|dH-mY|$, $m\geq 2$, is given by
$$\binom{d+n}{n} - 1-\sum_{i=0}^{m-1}
\binom{d+s-i}{d-i}
\binom{n-s-1+i}{i}.$$
In particular, the expected dimension of $|dH-2Y|$ is given by
\begin{equation}
\binom{d+n}{n} - 1-
\binom{d+s}{d}  -\binom{d+s-1}{d-1}(n-s)
\end{equation}

Consider now the system $\LL:=\LL_{n,d}(2^h)$, with $s+1\leq h$, and suppose 
that the first $s+1$ points span $Y=\PP^s$. Since $Y$ 
does not pass through all double points in $\LL$ we need to
study the system $|dH-2Y-\sum_{i=s+2}^h2P_i|$.  Then, for $s+1 \leq h$,  $Y=\PP^s$ is a 
$2-$special effect variety for $\LL_{n,d}(2^h)$ if
\begin{equation}\label{sineqsyslin2}
(s+1)(n+1)-\binom{d-s}{d}  -\binom{d-s-1}{d-1}(n-s)>0
\end{equation}
and
\begin{equation}\label{gensyslin}
\binom{d+n}{n} - 1-
\binom{d+s}{d}  -\binom{d+s-1}{d-1}(n-s)
-(h-s-1)(n+1)\geq 0.
\end{equation}

If we check the special inequality (\ref{sineqsyslin2}) for $d=2$ we obtain
$$(n+1)(s+1)-\frac{(s+1)(s+2)}{2}-(s+1)(n-s) > 0.$$
Simplifying we obtain $2n+2-s-2-2n+2s>0$ then $s>0$.

Let us consider the case $d\geq 3$. We write the special inequality as 
\begin{equation*}
\begin{split}
&(n+1)(s+1) - \binom{d+s}{s} - \binom{d-1+s}{d-1}(n-s)=\\
=&(n-s+s+1)(s+1)-\binom{d+s}{s} - \binom{d-1+s}{d-1}(n-s)= \\
=&(n-s)(s+1)+(s+1)^2-\binom{d+s}{s} - \binom{d-1+s}{d-1}(n-s).
\end{split}
\end{equation*}
Since $d\geq 3$ we have
$$\binom{d+s}{s} \geq \binom{s+3}{s} > (s+1)^2$$
and
$$\binom{d-1+s}{d-1} \geq\binom{s+2}{2} > s+1.$$
Thus, for $d\geq 3$, the special inequality (\ref{sineqsyslin2}) is false for every $s$ and $n$.

We pass now to study the equation (\ref{gensyslin}) assuming $d=2$. We obtain
\begin{equation}\label{s=h-1}
\frac{1}{2}n^2+\frac{3}{2}n+\frac{1}{2}s^2+\frac{1}{2}s-h(n+1)\geq 0.
\end{equation}
If we solve by respect to $s$ we find
$$s\geq \bigg\lfloor \frac{\sqrt{1-12n-4n^2+8hn+8h}}{2}-\frac{1}{2}\bigg\rfloor.$$
Thus we have the following
\begin{Prop}\label{syslingencase}
Let $\LL:=\LL_{n,d}(2^h)$ with $2\leq h \leq n$. Then $Y=\PP^s$ is a $2-$special effect variety for $\LL$ if $ \rho(n,h) \leq s \leq h-1$, where
$$\rho(n,h)=\begin{cases} 
\big\lfloor\frac{\sqrt{1-12n-4n^2+8hn+8h}}{2}-\frac{1}{2}\big\rfloor  & if \quad h> \frac{n^2+3n}{2(n+1)}\\
1 & otherwise
\end{cases}
$$
\end{Prop}

\begin{Cor}\label{syslin}
$Y=\PP^{h-1} \subset \PP^n$ is a $2-$special effect variety
for $\LL_{n,2}(2^h)$ for $2\leq h \leq n$ and $\forall n \geq 2$.
\end{Cor}

\begin{proof}
It follows easily from the proof of Proposition \ref{syslingencase} by
observing 
that formula (\ref{s=h-1}) is always verified for $s+1=h$, for $2\leq h \leq n$ and $\forall n \geq 2$.
\end{proof}


\subsection{Rational normal curves of degree $n$ in $\PP^n$}
Let $C_n \subset \PP^n$ be the image of $\PP^1$ under the $n-$Veronese embedding ($n>1$).
Once we have fixed the dimension $n$ of $\PP^n$, 
the virtual dimension of $|dH-mC_n|$ can be computed
using the generalized Riemann--Roch theorem on the blow-up of $\PP^n$ along $C_n$.
More generally, there are some classical results 
about the postulation of a multiple curve.
See, for example, the works of B. Segre (\cite{Segre}) and A. Franchetta (\cite{Franchetta}).

Since we are interested, for the moment,
only in the case $|dH-2C_n|$
we can use some interesting results given by A. Conca in \cite{Conca}.
Thus, for $d\geq 3$, one has
$$\nu(|dH-2C_n|)=\binom{d+n}{n} - 1-((d-1)n^2+2).$$
Supposing $h=n+3$ so that $C_n$ is fixed, the special inequality becomes
$$(n+1)(n+3)- (d-1)n^2-2>0.$$
If we expand the previous inequality we obtain 
\begin{equation*}
\begin{split}
(n+1)(n+3)-(d-1)n^2-2&=-(d-1)n^2-2+n^2+4n+3=\\
&=(2-d)n^2+4n+1>0.
\end{split} 
\end{equation*}
If we solve this equation with respect to $n$ we find 
$$\frac{2-\sqrt{2+d}}{d-2} < n < \frac{2+\sqrt{2+d}}{d-2} $$
then we restrict our solutions to $2 \leq  n < \frac{2+\sqrt{2+d}}{d-2}$. 
If we substitute the values of $d$ we find that the only possibilities are
$$\begin{array}{c|c|c}
d  & 3   & 4 \\ \hline
n  & 2,3,4& 2 
\end{array}$$
At this point,
we need to check $\nu(|dH-2C_n|)\geq0$.
Since 
$$\nu(|3H-2C_2|)=\nu(|3H-C_3|)=-1,$$
we exclude $d=3$ with $n=2,3$.
This concludes the proof of the following

\begin{Prop}\label{RNC}
Let $C_n\subset \PP^n$ be the rational normal curve, i.e. 
the image of $\PP^1$ under the $n-$Veronese embedding.
Then $C_n$ is a $2-$special effect variety for $\LL_{2,d}(2^{n+3})$
only when $(n,d)$ is $(2,4)$ or $(4,3)$. 
\end{Prop}

\begin{Rmk}\rm
It is easy to check that $C_n$ is not a $2-$special effect variety 
for $\LL:=\LL_{n,d}(2^h)$ if $h\not=n+3$.

Let us analyze first the case $h\leq n+2$. The conditions for 
the speciality of $C_n$ are
\begin{align}
\binom{d+n}{n} - 1-((d-1)n^2+2)\geq 0\label{1hn2} \\
h(n+1)-((d-1)n^2+2)>0 \label{2hn2}
\end{align}
Since $h\leq n+2$ and $d\geq 3$, in (\ref{2hn2}) we obtain
$$0<h(n+1)-((d-1)n^2+2)\leq (n+2)(n+1)-(2n^2+2)=3n-n^2$$
Thus the only possible value is $n=2$  and equation (\ref{2hn2}) becomes
$$2-4d+3h>0.$$
Since $C_2$ is the conic in $\PP^2$, we have to consider $d\geq 4$. 
Thus the previous equation has no solutions
for $h\leq 4$.

For $h\geq n+4$ the speciality inequality 
is the same as in the case of Proposition \ref{RNC}, 
hence the allowed values are
$$\begin{array}{c|c|c}
d  & 3   & 4 \\ \hline
n  & 2,3,4& 2 
\end{array}$$
But for this values, with the hypothesis $h\geq n+4$, one has
 $\nu(|dH-2C_n-\sum_{i=n+4}^h2P_i|)<0$.
\end{Rmk}

\subsection{The $\alpha-$special effect varieties and the Alexander--Hirschowitz Theorem}
The examples in the previous sections
fit with the Alexander--Hirschowitz Theorem.
In particular we can state the following

\begin{Thm}\label{asecandAH}
The Numerical Conjecture holds
for each of the special systems listed in Theorem \ref{P^nspecial}
\end{Thm}

\begin{proof}
It is enough to find a $\alpha-$special effect variety $Y$
for each of the special systems in the list of Thereom \ref{P^nspecial}.\\
The cases $\LL_{n,2}(2^h)$, $2\leq h \leq n$ follow
from Proposition \ref{syslin} considering $Y$ as the linear space $\PP^{s-1}$. \\
The conic in $\PP^2$ is a $2-$special effect curve
for $\LL_{2,4}(2^5)$ as shown in Example 3.8 in \cite{sev-1}.
It follows also from Proposition \ref{dH}, case $ii)$,
and from Proposition \ref{RNC}, case $(n,d)=(2,4)$. \\
The cases $\LL_{3,4}(2^9)$ and  $\LL_{4,4}(2^{14})$ are studied in Proposition \ref{dH}
and $Y$ is the quadric hypersurface respectively in $\PP^3$ and $\PP^4$. \\
Finally, using again Proposition \ref{RNC} case $(n,d)=(4,3)$,
we obtain  that the rational normal curve $C_4 \subset \PP^4$
is a $2-$special effect variety for  $\LL_{4,3}(2^7)$.
\end{proof}

\begin{Rmk}
Recently, A. Laface and L. Ugaglia proposed a conjecture for special
linear systems in $\PP^3$ (\cite{P3}).
Although an equivalence between this conjecture and the
Numerical Special Effect Conjecture is still unproved, it is easy to
see some interesting evidence. In fact, in the Laface--Ugaglia
Conjecture the speciality of a system $\LL$ in standard form (i.e. 
after performing a series of Cremona transformations) is related to the
existence of a quadric surface or a line in the base locus $Bs(\LL)$  which
makes the value of $\nu(\LL)$ lower. In other terms, both the
quadric or the line seem to be $\alpha-$special effect varieties.
\end{Rmk}


\section{$h^1-$Special effect varieties in $\PP^n$, $n\geq 3$}\label{h1high}
We turn now to analyzing $h^1-$special effect varieties in higher dimension.
Let $\LL:=\LL_{n,d}(-\sum_{i=1}^hm_iP_i)$ be a linear system of hypersurfaces
with general multiple base points
and let $X$ be the blow-up of $\PP^n$ at the points $\{P_i\}$.
Let $\tilde{\LL}$ be the strict transform of $\LL$.
In general, if confusion cannot arise, we will denote both $\LL$ and $\tilde{\LL}$ by $\LL$.
We recall that, if we denote by $\tilde{Y}$ the strict transform of a variety $Y\subset \PP^n$,
then we define $\LL-Y:=\LL\otimes\II_{\tilde{Y}}$.
The definition of the $h^1-$special effect variety
is slightly modified with respect to the planar case.

\begin{Def}\label{h1alta}
Let $\LL$ and $Y$ be as above with $Y$ irreducible.
Moreover, when codim$(Y,\PP^n)=1$,
we require $\OO_{\PP^n}(Y)\not\cong \LL$.
Then $Y\subset \PP^n$ is an $h^1-${\bf special effect variety}
for the system  $\LL$
if the following conditions are satisfied:
\begin{itemize}
\item[(a)] $h^0(\LL_{|Y})=0$;
\item[(b)] $h^0(\LL-Y)\not=0$;
\item[(c)] $h^1(\LL_{|Y})>h^2(\LL-Y)$.
\end{itemize}
\end{Def}

As in the planar case,
the speciality of the system $\LL$
follows from the previous conditions and from the standard exact sequence
\begin{equation}\label{standexaseq}
0 \to \LL-Y \to \LL \to \LL_{|Y} \to 0.
\end{equation}
In fact we have the following long exact sequence in cohomology:
\[
0 \to H^0(\LL-Y) \to H^0(\LL) \to H^0(\LL_{|Y})
\to  H^1(\LL-Y) \to H^1(\LL) \to H^1(\LL_{|Y}) \to \cdots
\]
Conditions $(a)$ and $(b)$
assure us that $H^0(\LL)\not= 0$,
while condition $(c)$ implies $H^1(\LL)\not= 0$.
Thus the existence of such $Y$
forces the system $\LL$ to have $h^0(\LL)\cdot h^1(\LL)\not= 0$
so that, by (\ref{chi}), $\LL$ is special.

\begin{Def}\label{cohomspec}
A special system arising from the existence of an $h^1-$special 
effect curve is called {\bf Cohomologically Special}.
\end{Def}

And again we can state a conjecture:

\begin{Conj}[(CSEC) ``Cohomological Special Effect'' Conjecture]\label{HCSEC}
A linear system  $\LL:=\LL_{n,d}(-\sum_{i=1}^hm_iP_i)$ 
with general multiple base points is special if and only if it is cohomologically special.
\end{Conj}

The $h^1-$special effect varieties seem easier to treat
than the $\alpha-$special effect varieties.
In fact we do not need to define the virtual dimension,
but we just work with elements in cohomology.
However, in several situations,
it is very difficult to compute some cohomology groups, in particular $h^2(\LL-Y)$. 

As in the case of $\alpha-$special effect varieties,
we do not have problems when $Y$ is a divisor
since $h^2(\LL-Y)=0$ if $\LL-Y$ is effective. Unluckily,
in this case, it can be difficult to study the behaviour of
$\LL_{|Y}$.

Instead, when codim$(Y,\PP^n)\geq 2$,
the groups $h^i(\LL-Y)$, $i=1,2$ can be computed on the blow-up of $\PP^n$ along $Y$, 
but we need a deep understanding of the geometry and cohomology of $Y$.  

We study now the situation in which $\LL:=\LL_{n,d}(2^h)$,
i.e. $\LL$ is a linear system with imposed double points.
The following Theorem is similar to Theorem \ref{asecandAH}.
Thus, also the Cohomological Special Effect Conjecture
fits with the Alexander--Hirschowitz Theorem.

\begin{Thm}\label{h1sevandAH}
The Cohomological Conjecture holds 
for each of the special systems listed in Theorem \ref{P^nspecial}.
\end{Thm}

\begin{proof}
We start with $\LL:=\LL_{n,2}(2^h)$ with $2\leq h \leq n$.
Let $Y$ be the span of the $h$ points $P_0, \dots P_{h-1}$ in the linear system $\LL$,
i.e. $Y=\PP^{h-1}$.
Since $\LL_{|Y}=\LL_{h-1,2}(2^h)$,
this system is clearly empty (see, for example, \cite{Miranda}) and
one has $h^0(\LL_{|Y})=0$
and the condition $(a)$ for $Y$ to be an $h^1-$special effect variety is satisfied.
Let $Z$ be the zero-dimensional scheme $\cup_{i=0}^{h-1}2P_i$ on $Y$;
then from the exact sequence
$$
\begin{array}{ccc}
0 \to & \II_Z(2) & \to \OO_{\PP^{h-1}}(2) \to \OO_Z \to 0\\
& \parallel \\
& \LL_{|Y}\\
\end{array}
$$
we obtain $h^1(\LL_{|Y})=\frac{h(h-1)}{2}$.

By Theorem \ref{P^nspecial} one has $h^0(\LL)\not=0$.
Since $h^0(\LL_{|Y})=0$ we conclude $h^0(\LL-Y)\not=0$,
then condition $(b)$ is satisfied.

Thus $Y$ will be an $h^1-$special effect variety 
if we prove $h^1(\LL_{|Y})>h^2(\LL-Y)$.
From the discussion, in \cite{Miranda},  about the matrices representing quadratic forms 
we easily compute 
$$h^1(\LL)=\frac{h(h-1)}{2}.$$
From the sequence (\ref{standexaseq}) we obtain
\begin{equation*}
0 \to H^1(\LL-Y) \to H^1(\LL) \stackrel{\theta}{\to} H^1(\LL_{|Y}) \to H^2(\LL-Y) \to 0.
\]
Since $h^1(\LL)=h^1(\LL_{|Y})$ it is enough to 
prove that $\theta$ is not the zero map.

We can suppose that the $h$ points
are the coordinate points 
$P_i=\lbrack 0,\dots, 1, \dots, 0\rbrack$, $i=0,\dots, h-1$
so that $Y$ has equation $x_h=\dots=x_n=0$.
Let $I_{Y}$ be the ideal of $Y$ in $\PP^n$
and let ${\bf m}_i,{\bf m}_{Y,i}$ be the ideals of $P_i$'s
respectively in $\PP^n$ and $Y$.
Let $I$ and $I'$ be respectively
the ideals $\cap_{i=0}^{h-1}{\bf m}_i^2$ and $\cap_{i=0}^{h-1}{\bf m}_{Y,i}^2$.\\
Since  ${\bf m}_{Y,i}=(x_0,\dots,\hat{x_i},\dots x_{h-1})$ for $i=0,\dots, h-1$,
the generators of $I'$ are
$$
\begin{cases}
x_kx_lx_m & k,l,m=0,\dots, h-1, k\not=l, k\not=m, l\not=m,\\
x_k^2x_l^2 & k,l=0,\dots, h-1, k\not=l.
\end{cases} 
$$
Moreover ${\bf m}_i=(x_0,\dots,\hat{x_i},\dots x_n)$  for $i=0,\dots, h-1$;
hence if $j=h,\dots,n$ then $x_j^2 \in \cap_{i=0}^{h-1}{\bf m}_i^2$.
Thus, after a straightforward computation we obtain
$$I=I_Y^2\cup (x_kx_lx_m: l,m=0,\dots, h-1, m\not=l, k=h, \dots, n ) \cup I'.$$
We denote by $\II_Y$, $\II$ and $\II'$ 
the ideal sheaves corresponding to the previous ideals.

Consider the following diagram:
\begin{equation}\label{DIA}
\begin{CD}
H^0(\OO_{\PP^n}(2)) @>\alpha>>H^0(\OO_{\PP^n}/\II(2)) = B \\
@VV{<x_h, \dots, x_n>}V  @VV{<x_h, \dots, x_n>}V\\
H^0(\OO_Y(2)) @>{\alpha_Y}>>H^0(\OO_Y/\II'(2)) = B_Y
\end{CD}
\end{equation}
We call $\sigma$ the map $B \to B_Y$ in the previous diagram
(the map on the right-side). 
This map is given by the equations of $Y$.

From the previous computations of the ideals 
$I$ and $I'$ we know that $\sigma$ is surjective.
Since $H^1(\LL_{|Y})\not= \emptyset$,
there exists an $\eta \in B_Y$ such that  
$\eta \notin \mbox{Im}\alpha_Y$. Let $\lbrack \eta\rbrack$ 
be the image of $\eta$ in $H^1(\LL_{|Y})$.
By the surjectivity of $\sigma$,
$\eta$ comes from an element $\eta_0 \in B$.
Since diagram (\ref{DIA}) is commutative,
$\eta_0$ does not lies in $\mbox{Im}{\alpha}$,
because otherwise $\eta \in \mbox{Im}\alpha_Y$.
Thus we conclude that $\theta$ sends $\lbrack \eta_0\rbrack \in H^1(\LL)$
to $\lbrack \eta \rbrack \in H^1(\LL_{|Y})$, i.e. $\theta\not\equiv 0$.
Hence $h^1(\LL_{|Y})>h^2(\LL-Y)$. 

We turn to analyzing the cases $\LL_{2,4}(2^5)$, $\LL_{3,4}(2^9)$ and $\LL_{4,4}(2^{14})$.
We can treat them in an unified way just writing $\LL:=\LL_{n,4}(2^s)$,
where $s=\left(\begin{smallmatrix} 2+n \\ n \end{smallmatrix}\right)-1$ and $n=2,3,4$.
Let $Y$ be the divisor corresponding to $\LL_{n,2}(1^s)$,
i.e the conic in $\PP^2$ through $5$ points
and the quadric in $\PP^3$ and $\PP^4$ respectively through $9$ and $14$ points.
Since $\LL-Y=\LL_{n,2}(1^s)$ we have
$$h^0(\LL-Y)=1,$$
$$h^i(\LL-Y)=0 \mbox{ for } i\geq 1.$$ 
By Theorem \ref{P^nspecial} we know that $h^0(\LL)=1$ and 
$$h^1(\LL)=\begin{cases}
1 & \mbox{ for } n=2,4\\
2 & \mbox{ for } n=3 
\end{cases}
$$
Finally, by sequence (\ref{standexaseq}) we conclude $h^0(\LL_{|Y})=0$
and $h^1(\LL_{|Y})=h^1(\LL)>0$.
Hence the given $Y=\LL_{2,n}(1^s)$
is an $h^1-$special effect variety for the system  $\LL:=\LL_{n,4}(2^s)$, $n=2,3,4$.

The last case to treat is $\LL_{4,3}(2^7)$.
Let $Y$ be the rational normal curve of degree $4$
passing through the seven double points described in \cite{CiHi}.
Since $\LL\cdot Y=-2$ we have
$$h^0(\LL_{|Y})=0 \qquad h^1(\LL_{|Y})=1.$$
Moreover, by Theorem \ref{P^nspecial}, we know that 
$$h^0(\LL)=1 \qquad h^1(\LL)=1.$$
Thus we obtain $h^0(\LL-Y)=1$.
To conclude the proof we need to show $h^2(\LL-Y)=0$. By the sequence 
\begin{equation*}
0 \to H^1(\LL-Y) \to H^1(\LL) \stackrel{\theta}{\to} H^1(\LL_{|Y}) \to H^2(\LL-Y) \to 0,
\]
it is enough to prove that $\theta$ is surjective. 
If we tensor by $\OO_{\tilde{Y}}$ the sequence (\ref{standexaseq}) we obtain
$$
\begin{array}{cccccccccc}
& 0 & \to &\LL & \to &\OO_{\bPP^4}(3) & \to & \OO_{\sum 2E_i} & \to & 0\\
& & & \downarrow & & \downarrow & & \downarrow \\
0 \to &\mbox{Tor}_1(\OO_{\sum2E_i},\OO_{\tilde{Y}}) & \to &\LL\otimes \OO_{\tilde{Y}} & \to &\OO_{\tilde{Y}}(3) & \to & \OO_{\sum 2Q_i} & \to & 0\\
\end{array}
$$
where $Q_i=\tilde{Y}\cap E_i$, $i=1,\dots,7$.
Since $\LL\otimes \OO_{\tilde{Y}}$ corresponds to an invertible sheaf on $\PP^1$,
it cannot have torsion,
thus $\mbox{Tor}_1(\OO_{\sum2E_i},\OO_{\tilde{Y}})=0$ and we write 
\begin{equation}\label{commdig}
\begin{array}{ccccccccc}
0 & \to &\LL & \to &\OO_{\bPP^4}(3) & \to & \OO_{\sum 2E_i} & \to & 0\\
& & \downarrow & & \downarrow & & \downarrow \\
0 & \to &\LL\otimes \OO_{\tilde{Y}} & \to &\OO_{\tilde{Y}}(3) & \to & \OO_{\sum 2Q_i} & \to & 0.\\
& & \parallel \\
& & \LL_{|Y}
\end{array}
\end{equation}
Since every diagram in (\ref{commdig}) is commutative,
 when we pass to cohomology we obtain the following commutative diagram
\[
\begin{CD}
H^0(\OO_{\sum 2E_1}) @>{\delta_1}>> H^1(\LL)\\
@V{\alpha}VV @V{\theta}VV\\
H^0(\OO_{\sum 2Q_1}) @>{\delta_{1,Y}}>> H^1(\LL_{|Y})\\
\end{CD}
\]
where $\delta_1$ and $\delta_{1,Y}$ are the connection homomorphisms.
We can observe that $\delta_1$ and $\delta_{1,Y}$ are surjective
because $H^1(\OO_{\bPP^4}(3))$ and $H^1(\OO_{\tilde{Y}}(3))$ are zero.
Moreover $\alpha$ is surjective too.
As a matter of fact, let $f_i\in k\lbrack\dots, x_j ,\dots \rbrack$ be the polynomial defining $E_i$.
We can fix our attention on a single polynomial $f_0$.
Let $I(\tilde{Y})$ be the ideal of $\tilde{Y}$,
thus $k \lbrack \dots, x_j, \dots \rbrack/I(\tilde{Y}) =k \lbrack t \rbrack$.
Since $f_0$ is not tangent to $\tilde{Y}$ we have 
$$f_0 \mod I(\tilde{Y}) = t+o(t^2).$$
Thus the map
$$k \lbrack\dots, x_j,\dots \rbrack/(f_0^2) \to k \lbrack t \rbrack/(f_0^2)$$
is surjective. Finally, from the surjectivity of $\delta_{1,Y}\circ\alpha=\theta\circ\delta_1$
it follows that $\theta$ is surjective.
\end{proof}

\begin{Rmk}
From Theorems \ref{asecandAH} and \ref{h1sevandAH} we can notice that
 each $\alpha-$special effect variety for special systems in Theorem \ref{P^nspecial}
 is an $h^1-$special effect variety too for the same system.
However in $\PP^n$, $n\geq 3$, this is not true in general,
as will be shown in Example \ref{alphasih1no}. 
\end{Rmk}


\section{More examples of special effect varieties in $\PP^n$}\label{sevexs}

We collect in this section some special systems arising from the existence of different
kind of special effect varieties. In particular we show a variety for
the Laface--Ugaglia example (\cite{LafUga}) which is both $\alpha-$special effect and
$h^1-$special effect.

\begin{Ex}\rm {\bf (Homogeneous special systems in $\PP^n$)}
Let $Y$ be a linear space $\PP^s \subset \PP^n$.
It is not difficult to construct a family of homogeneous special systems $\LL_{n,d}(m^{s+1})$
with $Y$ as a special effect variety.
Again we underline that the study of special effect varieties
can turn in a pure combinatorial problem.

As an example we just consider $Y=\PP^1\subset \PP^3$,
i.e. $s=1$ and $n=3$. 
In this case, we write the special inequality as
\begin{equation*}
\begin{split}
&2\binom{m+2}{3} -  \sum_{i=0}^{m-1}\binom{d+1-i}{d-i}\binom{i+1}{i}=\\
=& \frac{m(m+1)(m+2)}{3}- \sum_{i=0}^{m-1}(d+1-i)(i+1)=\\
=& \frac{m(m+1)(m+2)}{3}- \sum_{i=0}^{m-1}\left[(d+1)+di-i^2)\right]=
\end{split}
\end{equation*}
\begin{equation*}
\begin{split}
=& \frac{m(m+1)(m+2)}{3} - m(d+1)-\frac{d(m-1)m}{2}+\frac{(m-1)m(2m-1)}{6}= \\
= & \frac{m}{6}(4m^2+3m-1-3d-3md)>0
\end{split}
\end{equation*}
Thus we ask for $4m^2+3m-1-3d-3md>0$ and we obtain that
$\PP^1$ is an $m-$special effect variety for $\LL_{3,d}(m^2)$ if $m\leq d<\frac{4m-1}{3}$.
In a similar way we can prove that $Y=\PP^2$ is an $m-$special effect
variety for $\LL_{3,d}(m^3)$ if 
$$m \leq d \leq \frac{m}{2}-2+\frac{\sqrt{84+108m+33m^2}}{6}.$$
\end{Ex}

\begin{Ex}\rm {\bf (Rational curves in $\PP^3$)} 
Let $Y$ be a smooth rational curve in $\PP^3$ and define $X$ as the blow-up of $\PP^3$ along $Y$. Thus we have the diagram
\begin{equation}\label{blowup}
\begin{CD}
R @>j>>X \\
@VV{g}V  @VV{\pi}V\\
Y @ >i>>\PP^3
\end{CD}
\end{equation}
where $R=\PP({\mathcal{N}}_{C|\PP^3})$ is the exceptional divisor along $Y$.
Let $\tilde{H}$ be the pull-back via $\pi$ of a general hyperplane
section of $\PP^3$.

The virtual dimension of $|dH-2Y|$ can be computed as $\chi(\OO_X(d\tilde{H}-2R))$ on $X$.
Using the generalized Riemann--Roch theorem (\cite{Fulton} pages 286--295) we obtain
$$\chi(\OO_X(D))=\frac{1}{12}D\cdot (D-K)\cdot(2D-K)+\frac{1}{12}D\cdot c_2+1$$
where $K:=K_X$. 
Since $c_1(X)=\pi^*(c_1(\PP^3))-R$ (for the proof, see \cite{GrifHar},
page 608) we have $K=-4\tilde{H}+R$.

Suppose that $Y$ has degree $e$. Since $c_2(\PP^3)=6H^2$,
from \cite{GrifHar} (Lemma at pages 609--610), we obtain
$c_2=(6+e)\tilde{H}^2-4\tilde{H}\cdot R$.
Thus we can write $\chi(\OO_X(d\tilde{H}-2R))$ as
\begin{equation*}
\begin{split}
\chi(\OO_X(d\tilde{H}-2R))&=\frac{1}{12}(d\tilde{H}-2R)\cdot((d+4)\tilde{H}-3R)\cdot((2d+4)\tilde{H}-5R)+\\
&+\frac{1}{12}(d\tilde{H}-2R)\cdot \left[(6+e)\tilde{H}^2-4\tilde{H}\cdot R\right]+1
\end{split}
\end{equation*}
We recall that $\tilde{H}^3=1$, $\tilde{H}^2 \cdot R =0$,
$\tilde{H}\cdot R^2 = (\tilde{H} \cdot R) \cdot R = eF \cdot R = -e$
and $R^3=2-4e$ (the last one can be computed by using Proposition at
page 606 in \cite{GrifHar}).
Using these results we obtain
$$\chi(\OO_X(d\tilde{H}-2R))=\frac{1}{6}d^3-5+d^2-3de+\frac{11}{6}d+4e=\binom{d+3}{3} -3de+4e-5.$$

A rational curve of degree $e$ in $\PP^3$ can be defined by four
polynomials of degree $e$. 
The set of their coefficients defines a projective space of 
dimension $4(e+1)-1-\mbox{Aut}(\PP^1)=4e+4-1-3=4e$. 
Since a simple point imposes two conditions on the equations of the
curve, we can use, 
in a first analysis, the bound $4e\geq 2h$ where $h$ is the number of
points. 
Obviously we need to check if the $h$ points are in general position. 
In fact this is not a consequence of the previous bound.

Thus we want $Y$ passing through the $h$ points of $\LL_{3,d}(2^h)$.
The conditions for the $2-$speciality of $Y$ are
\begin{align}
&2e\geq h,\label{g1}\\
&\binom{d+3}{3} -3de+4e-6\geq0,\label{g2}\\
&4h> 3de-4e+5.\label{g3}
\end{align}
From (\ref{g1}) and (\ref{g3}) we obtain $8e\geq 4h>3de-4e+5$, then $(12-3d)e>5$. 
This forces $d\leq 3$ and we finally find
\begin{itemize}
\item[(a)] the line is a $2-$special effect curve for $\LL_{3,2}(2^2)$,
\item[(b)] the conic is a $2-$special effect curve for $\LL_{3,2}(2^3)$.
\end{itemize}
Case $(a)$ was already discovered in Proposition \ref{syslin}. 
Case $(b)$ exhibits a new $2-$special effect variety for the system $\LL_{3,2}(2^3)$:
the other one was the plane $\PP^2$, by Proposition \ref{syslin}. 
Moreover, since we fix only three points, the conic can move and it
fills exactly a $\PP^2$. 
More generally it is possible to prove that the special system $\LL_{n,2}(2^h)$, 
for a fixed $h$, has at least two special effect varieties: 
the linear space $\PP^{h-1}$ and the rational normal curve $C_{h-1}\subset \PP^{h-1}$. 
\end{Ex}

\begin{Ex}\rm {\bf (Particular unions of lines in $\PP^n$ )} Let $Y\subset
  \PP^n$ be
  the union of the $\left(\begin{smallmatrix} n+1 \\ 2
  \end{smallmatrix}\right)$ lines passing through the $n+1$ coordinate points 
$P_i=[0,0,\dots, 1, \dots 0]$, for $i=0, \dots, n$. 
Then $Y$  is a $(n-1, \dots, n-1)-$special effect configuration for 
$\LL_{n,n+1}(n^{n+1})$, $n\geq 3$ (the proof is left to the reader).
\end{Ex}

\begin{Ex}{\bf (The Laface--Ugaglia Example)} 
In \cite{LafUga} Laface and Ugaglia show a counterexample to a
conjecture presented in \cite{Ciliberto} which requires, for a special system, the
existence of a rational curve in the base locus.
Laface and Ugaglia analyzed the linear system $\LL:=\LL_{3,9}(-6P_0-\sum_{i=1}^84P_1)$. 
It splits as $Q+\LL'$, where $\LL'=\LL_{3,7}(-5P_0-\sum_{i=1}^83P_i)$ and $Q$ is the quadric 
in $\PP^3$ passing through $P_0, \dots, P_8$. 
Then $\LL$ is special because $\nu(\LL)=3$ while $\nu(\LL')=4$. 
After that, they proved that the only curve contained in the base locus of $\LL$ 
is a curve $C\subset Q$ of genus $2$ given by the intersection of $Q$ with the generic 
element in $\LL'$.

From the previous considerations we see that $Q$ is a $1-$special effect variety for $\LL$. 
As a matter of fact we have
\begin{itemize}
\item[(i)] $\nu(|Q|)=0$,
\item[(ii)] $\nu(\LL-Q)=4>3=\nu(\LL)$.
\item[(iii)] $\nu(\LL-Q)=\nu(\LL_{3,7}(-5P_0-\sum_{i=1}^83P_i))=4$,
\end{itemize}
Consider now the restricted system $\LL_{|Q}=|9L_1+9L_2-6P_0-\sum_{i=1}^84P_1|$, 
where $L_1,L_2$ are the generators of Pic$(Q)$. 
$\LL_{|Q}$ is empty of virtual dimension $-2$ (see the Appendix in \cite{LafUga} for the proof). 
Hence $h^0(\LL_{|Q})=0$ and working on the blow-up of $Q$ at the $P_i$'s we obtain
\vskip0.2cm
\noindent
$h^1(\LL_{|Q})=$\vskip0.1cm
$\quad=h^2(\LL_{|Q})-\frac{(9\tilde{L}_1+9\tilde{L}_2-6E_0-\sum_{i=1}^84E_i)(11\tilde{L}_1+11\tilde{L}_2-7E_0-\sum_{i=1}^85E_i)}{2}-1=$\vskip0.1cm
$\quad=h^2(\LL_{|Q})+2-1\geq 1.$
\vskip0.2cm
Finally $h^0(\LL-Q)=h^0(\LL_{3,7}(-5P_0-\sum_{i=1}^83P_i))=4$  
and $h^2(\LL-Q)=0$. Then $Q$ verifies the conditions to be an $h^1-$special effect variety too.
\end{Ex}

\begin{Ex}{\bf (A $1-$special effect variety that is not an
    $h^1-$special effect variety)}\label{alphasih1no} 
Consider the system $\LL:=\LL_{3,d}(m^3)$ and let $Y\subset \PP^3$ be the plane through the three points in $\LL$.
Writing down the conditions of speciality we see that $Y$ is an
$1-$special effect variety for $\LL$ if
\[ 
\begin{cases}
\frac{-3+\sqrt{1+12m^2+12m}}{2}>d \\
 d^3+3d^2+2d-6\geq 3m^3-3m
\end{cases}
\]
Consider now the system $\LL:=\LL_{3,6}(4^3)$. 
For the previous computation, $Y=\PP^2$ is a $1-$special effect variety for 
$\LL$ and $\LL$ is special.
One has $\nu(\LL)=23$ and $\nu(\LL-Y)=25$ as we expect by the special effect of $Y$. 
Moreover we can observe that $Y$ is not a $2-$special effect variety
since $\nu(\LL-2Y)=22$. 

If we restrict the system to $Y$ we obtain the planar system $\LL_{|Y}=\LL_{2,6}(4^3)$. 
This system is special, so both $h^0(\LL_{|Y})$ and  $h^1(\LL_{|Y})$ are different from zero. 
Hence $Y=\PP^2$ does not satisfy condition $(a)$ to be an
$h^1-$special effect variety for $\LL$.

However there is an $h^1-$special effect variety for the system $\LL$.
As a matter of fact, if we compute the effective dimension of the system $\LL$ by a computer 
algebra program (e.g Maple) we discover $\dim(\LL)=26$ then $\LL-Y$ represents a subsystem of the system 
of divisors of degree $6$ with three points of multiplicity $4$ (i.e $\LL$ does 
not split as $Y+\LL_{3,5}(3^3)$). 
Hence the generic element $D\in \LL$ cannot be written 
as the sum of $Y$ and of elements in $\LL_{3,5}(3^3)$. Suppose the points $P_i$'s are 
$P_1:=[0,1,0,0]$, $P_2:=[0,0,1,0]$ and $P_3:=[0,0,0,1]$ and the coordinates are 
$x_0, \dots x_3$. Then $Y$ is the plane defined by $x_0=0$. Moreover
$\LL$ is generated by
the span of
$<Y+\LL_{3,5}(3^3),F>$ where $F$ 
is $(x_1x_2x_3)^2$ (as we expect from Theorem 2.4 in \cite{CiMi2}).

We can observe that $F$ contains twice the lines $L_{ij}:=\overline{P_iP_j}$, $i,j=1,2,3$ 
and $i \not= j$. Moreover, every element in $Y+\LL_{3,5}(3^3)$ contains the same 
lines with multiplicity at least $2$. Thus $Y'=L_{12}+L_{13}+L_{23}$ is a $(2,2,2)-$special effect configuration for $\LL$.
Finally it is easy to check that if we just consider one of the previous 
lines $L_{ij}$ we obtain that $L_{ij}$ is an $h^1-$special effect
variety for both $\LL$ and $\LL_{|Y}$.
\end{Ex}


\section{$\alpha-$Special effect varieties in the product of projective spaces}\label{alphacross}

We show now several examples of $\alpha-$special effect varieties
on $X=\PP^{n_1}\times \cdots \times \PP^{n_t}$
with $t\geq 2$ and $n_i\geq 1$ for $i=1,\dots, t$.
We treat only the case case $m=\alpha=2$
and we suppose that the special effect variety $Y$ is a divisor on $X$.
Surely this does not exhaust all possible special effect varieties 
(and special linear systems) on $X$,
but we will observe at the end of the Section how our results
fit with the ones of Catalisano, Geramita and Gimigliano
on secant varieties of products of projective spaces
(\cite{CaGeGi2}, \cite{CaGeGi3}, \cite{CaGeGi4}).

\vskip0.5cm
{\bf Notation.}
Let $r$ be a positive 
integer. For any integer $z$ we define 
$(r)_{(z)}$ as follows
$$(r)_{(z)}:=
\begin{cases}
\Pi_{i=1}^z(r+i) & \mbox{ if } z>0\\
1 & \mbox{ if } z=0\\
0 & \mbox{ if } z<0
\end{cases}
$$
We have the following fact: let $r,s$ and $t$ be positive
integers, one has the equality
\begin{equation}\label{conto1}
(r+s)_{(t)}=(s)_{(t)}+r\left(\sum_{i=1}^t(s)_{(i-1)}(r+s+i)_{(t-i)}\right).
\end{equation}
Since each term $(s)_{(i-1)}(r+s+i)_{(t-i)}$ in the summation is greater than $(s)_{(t-1)}$ we can write
the following inequality
\begin{equation}\label{conto2}
(r+s)_{(t)}\geq (s)_{(t)}+rt\left((s)_{(t-1)}\right)= (s)_{(t-1)}(s+t+rt)
\end{equation}
\vskip0.5cm

Let $Y$ be a divisor of multidegree $(e_1, \dots. e_t)$ on $\PP^{n_1}\times \cdots \times \PP^{n_t}$ and consider
the system $\LL:=\LL_{(d_1, \dots, d_t)}(2^h)$ of divisor of multidegree 
$(d_1, \dots, d_t)$ passing through $h$ general double points. We
require that $Y$ passes through the $h$ points of $\LL$.
Then $Y$ is a $2-$special effect varieties for $\LL$ if
\begin{align}
&\Pi_{i=1}^t\binom{e_i+n_i}{n_i}
-1 \geq h;\label{prodcon1}\\
&\Pi_{i=1}^t\binom{d_i-2e_i+n_i}{n_i}
> \Pi_{i=1}^t \binom{d_i+n_i}{n_i}
-h(\sum_{i=1}^tn_i+1)\label{prodcon3}\\
&\Pi_{i=1}^t\binom{d_i-2e_i+n_i}{n_i}
\geq 1, \mbox{ i.e. }
d_i \geq 2e_i, \mbox{ for } i=1, \dots, t\label{prodcon2}.
\end{align}

We apply the same argument of Proposition \ref{dH}.
Again, the previous conditions give us the bounds on the number of points $h$:
\[\label{boundcross}
 \Pi_{i=1}^t\binom{e_i+n_i}{n_i}-1 \geq h >  \frac{1}{\sum_{i=1}^t n_i+1}\left[ \Pi_{i=1}^t\binom{d_i+n_i}{n_i}-\Pi_{i=i}^t\binom{d_i-2e_i+n_i}{n_i}\right]
\]
\noindent
so that we can study when the function
\begin{multline}\label{generalcross}
\varphi(d_1,\dots, d_t, e_1, \dots, e_t, n_1, \dots, n_t):=\\
=\Pi_{i=1}^{t}\binom{d_i+n_i}{n_i} - \Pi_{i=1}^{t}\binom{d_i-2e_i+n_i}{n_i} -  \left( \Pi_{i=1}^{t}\binom{e_i+n_1}{n_i} -1\right) \left( \sum_{i=1}^{t}n_i+1\right)
\end{multline}
\noindent
is negative.

\begin{Lem}[Numerical Lemma]\label{NL}  
Let  $\varphi(d_1,\dots, d_t, e_1, \dots, e_t, n_1, \dots, n_t)$ be
defined as in (\ref{generalcross}). Then the function  
\[
\eta(e_1, \dots, e_t, n_1, \dots, n_t):=\varphi(2e_1,\dots, 2e_t, e_1,
\dots, e_t, n_1, \dots, n_t)
\]
is non-decreasing in the $n_i$'s for
\begin{itemize}
\item[a)] $t=2$ with $e_i \geq 2$ and $n_i\geq 2$, $i=1,2$.
\item[b)] $t\geq 3$ with $e_i\geq 1$ and $n_i \geq 1$, $i=1, \dots, t$;
\end{itemize}
\end{Lem}

\begin{proof}
By definition of $\varphi$, we have
\begin{equation*}
\begin{split}
\eta&:=\varphi(2e_1, \dots, 2e_t, e_1, \dots, e_t, n_1, \dots, n_t)=\\
&=\Pi_{i=1}^{t}\binom{2e_i+n_i}{n_i} - \left( \Pi_{i=1}^{t}\binom{e_i+n_1}{n_i} -1\right) \left( \sum_{i=1}^{t}n_i+1\right)-1.
\end{split}
\end{equation*}
For the simmetry of $\eta$, it is enough to prove the lemma for one
$n_{i_o}$, with $i_o \in \{1, \dots, t\}$.

$\boxed{\mbox{case $(a)$}}$ We  prove that $\varphi(2e_1,2e_2,e_1,e_2,n_1,n_2)$ 
is increasing in $n_1$.
Thus, after a tedious computation, one has
\begin{align*}
&\gamma_{n_1}:=\varphi(2e_1,2e_2,e_1,e_2,n_1+1,n_2)-\varphi(2e_1,2e_2,e_1,e_2,n_1,n_2)=\\
&=\mbox{{\small{$P \cdot${\huge{$\lbrack$}}$(n_1+e_1)_{(e_1)}(n_2+e_2)_{(e_2)}-$}}} \\
&- \mbox{{\small{$(e_1)_{(e_1-1)}(e_2)_{(e_2)}\left((n_1+e_1+1)(n_1+n_2+2)-(n_1+1)(n_1+n_2+1)\right)${\huge{$\rbrack$}}}}}+1> \\
&> P \cdot \left((e_1)_{(e_1-1)}(e_2)_{(e_2)}\right)\left[e_1e_2n_1n_2-2e_2n_1-2e_2 \right]+1\geq0
\end{align*}
\noindent
where
$$P:=\frac{(n_1+1)_{(e_1-1))}(n_2)_{(e_2)}}{(2e_2)!(2e_1-1)!}$$
and, for the inequality we use (\ref{conto2}) with $r=n_i$, $s=t=e_i$, $i=1,2$.

Thus $\varphi$ is increasing in $n_i$ and $(a)$ follows.

$\boxed{\mbox{case $(b)$}}$ As in case $(a)$ we look for a good way to collect terms in 
\[
\gamma_{n_{i_0}}:= \eta(2e_1,\dots,2e_t, \dots, n_{i_0}+1,\dots, n_t)-\eta(2e_1,\dots,2e_t, \dots, n_{i_0},\dots, n_t)
\]

Using again (\ref{conto2}) one has
\begin{equation}
\begin{split}
\gamma_{n_{i_0}}>& \tilde{P} \cdot \left[ \frac{2e_{i_0}\Pi_{s=1}^t(n_s+e_s)_{(e_s)}}{\left(\Pi_{s=1}^t(e_s)_{(e_s)}\right)(n_{i_0}+e_{i_0}+1)}-\left(\sum_{i=1}^{t}n_i+2\right)\right]+1 =\\
= & \tilde{P}\cdot \left[C(e_i) \right]+1
\end{split}
\end{equation}
where
\[
P:=\frac{\Pi_{s=1}^t(n_s)_{(e_s)}}{(n_{i_0}+1)(\Pi_{k=1}^{t}(e_k!))}
\quad \mbox{ and } \quad \tilde{P}:= (n_{i_0}+e_{i_0}+1)\cdot P.
\]
Now we can apply the same argument of claim of Lemma \ref{LN1} and we
obtain that the term $C(e_i)$ is increasing in $e_i$ and
$\gamma_{n_{i_0}}$ too. 

If we substitute $e_1=\dots= e_t=1$ in $\gamma_{n_{i_0}}$ we obtain
\begin{equation}\label{con2et}
\gamma_{n_{i_0}}^* = \frac{P(e_i=1)}{2^{t-1}}\cdot \left[\sum_{m=1}^t\left( 2^{t-m}\sum_{|I|=m}n_{|I|}  \right)-2^{t-1}n_{i_0}-2^{t-1}\right]+1
\end{equation}
where $n_{|I|}=n_{i_1}\cdots n_{i_m}$ if $I=\{i_1, \dots, i_m \}$.

If $t\geq4$ we have at least two terms of the form $2^{t-2}n_jn_{i_0}$ and four terms of the form $2^{t-3}n_jn_kn_l$. Then we have
$$2^{t-2}\sum n_jn_{i_0}\geq 2^{t-1}n_{i_0}$$
and
$$2^{t-3}\sum n_jn_kn_l\geq 2^{t-1}.$$
Thus the expression between square brackets in (\ref{con2et}) is
always positive and then $\eta$ is increasing on $n_{i_0}$.

When $t=3$, we obtain, for example for $i_0=1$,
$$\gamma_{n_1}^*=\frac{P}{4}\left[n_1n_2n_3+2n_1n_2+2n_1n_3+2n_2n_3-4n_1-4\right]+1$$
\noindent
and the expression between square brackets is positive except for $n_1=n_2=n_3=1$, but, for these values we have
$$\gamma_{n_1}=\frac{P}{4}\left[n_1n_2n_3+2n_1n_2+2n_1n_3+2n_2n_3-4n_1-4\right]_{|n_i=1}+1=\frac{4}{4}\left[-1\right]+1=0$$
\noindent
then, also for the case $t=3$, $\eta$ is non-decreasing in $n_i$, $i=1,2,3$. Hence $(b)$ is proved.
\end{proof}

\begin{Prop}\label{prodprop}
Let  $\varphi(d_1,\dots, d_t, e_1, \dots, e_t, n_1, \dots, n_t)$ be
defined as in (\ref{generalcross}). Then
\begin{itemize}
\item[(a)] If $t=2$ then $\varphi\geq 0$ for $e_1,e_2\geq2 $, 
  for $n_1,n_2\geq 2$ and $d_i\geq 2e_i$, $i=1,2$.
\item[(b)] If $t=3$ then  $\varphi(d_1, d_2, d_3, e_1, e_2, e_3, n_1,
  n_2, n_3)\geq0$, 
except for 
$(d_1,d_2,d_3)=(2,2,2)$, $(e_1,e_2,e_3)=(1,1,1)$ and
$(n_1,n_2,n_3)=(1,1,\gamma)$ with $\gamma \leq 3$
\item[(c)] If $t\geq 4$ then $\varphi(d_1,\dots, d_t, e_1, \dots, e_t,
  n_1, \dots, n_t)\geq0$, 
for $e_i,n_i\geq 1$ and $d_i\geq 2e_i$, $i=1, \dots t$;
\end{itemize}
\end{Prop}

\begin{proof}
Since the function $\varphi$ is non-decreasing in $d_i$ we can start from the value $d_i=2e_i$, $i=1,\dots,t$.
Then, using the previous lemma it is enough to substitute the minimal values of $n_i$, $i=1, \dots, t$ in $\eta$ and then 
study the positivity of this easier function.

As an example, we prove $(c)$. In this case, after the
substitution of $n_i=1$ in $\eta$ we obtain
$$\varphi=\Pi_{i=1}^t(2e_i+1)-(t+1)\Pi_{i=1}^t(e_i+1)+t.$$
The previous expression is increasing in $e_i$, $i=1, \dots t$ and, if we finally substitute $e_1=\dots =e_t=1$ we obtain
$$\varphi=3^t-2^t(t+1)+t$$
and it is positive for $t\geq4$. Hence $(c)$ is proved.
\end{proof}

We now search for $2-$special effect divisors on $\PP^{n_1} \times
\cdots \times \PP^{n_t}$. We analyze first the case $t=2$.

\begin{Prop}\label{listpapb}
Let $\LL:=\LL_{(d_1,d_2)}(2^h)$ be a linear system of bidegree $(d_1,d_2)$
on $X=\PP^{n_1}\times \PP^{n_2}$
passing through $h$ double points in general position, with $d_1\cdot d_2\not=0$.
Let $Y \subset \PP^{n_1}\times \PP^{n_2}$ be a divisor of bidegree $(e_1,e_2)$,
with $e_i\not= 0$ for at least one $i$. Moreover we require that $Y$ passes simply
through the $h$ points in $\LL$.
Then $Y$ is a $2-$special effect variety for $\LL_{(d_1,d_2)}(2^h)$
 in the following cases
\[
\begin{array}{lccc}
\hline
\PP^{n_1}\times\PP^{n_ 2} & (d_1,d_2) & (e_1,e_2) & h\\\hline
\PP^{1}\times \PP^{1} & (2,2e_2) &  (1,e_2)  & 2e_2+1 \\
\PP^{1}\times \PP^{1} & (2e_1,2) &  (e_1,1) & 2e_1+1\\
\PP^{1}\times \PP^{n_2} & (2e_1,2) & (e_1,1) & m_1(e_1,n_2)\leq h\leq M_1(e_1,n_2) \\
\PP^{2}\times \PP^{n_2}& (2,2) & (1,1) & m_2(n_2)\leq h\leq M_2(n_2) \\
\PP^{3}\times \PP^{3} &  (2,2) & (1,1) & 15 \\
\PP^{3}\times \PP^{4}& (2,2) & (1,1) & 19 \\
\hline
\end{array}
\]
where
$$
\begin{array}{ll}
m_1(e_1,n_2):= \lfloor \frac{(2e_1+1)(n_2+1)}{2} \rfloor  & m_2(n_2):=\lfloor \frac{3n_2^2+9n_2+5}{n_2+3} \rfloor \\ 
 & \\
M_1(e_1,n_2):=e_1n_2+e_1+n_2  & M_2(n_2):= 3n_2+2.
\end{array}
$$
\end{Prop}

\begin{proof}
We start with the case $e_1\cdot e_2\not=0$. From Proposition
\ref{prodprop} we can restrict our analysis to $e_1+e_2\leq 3$ or
when at least one between $n_1$ and $n_2$ is equal to one.
We divide the proof in three steps,
analyzing some different situations.
\vskip0.2cm
{\bf Step 1:} $n_1=n_2=1$.
\begin{equation}\label{casor11}
\varphi(2e_1,2e_2,e_1,e_2,1,1)= e_1e_2-e_1-e_2
\end{equation}
\noindent
Obviously, (\ref{casor11}) is negative only  for 
$$\begin{cases}e_1=1 \\ e_2\geq 1 \end{cases}  \mbox{ or } \begin{cases} e_1\geq 1 \\ e_2=1
\end{cases}$$
\noindent
For the symmetry of the function, we can fix our attention on the case $e_1=1$. 
The only possible value for $(d_1,d_2)$ is $(2,2e_2)$; 
as a matter of fact, if $(d_1,d_2)=(2+i,2e_2+j)$ with $i,j\geq 1$,  we have
\[
\varphi(2+i,2e_2+j,1,e_2,1,1)=2e_2i+2j-1 > 0 \qquad \forall i,j\geq1.
\]
\noindent 
Thus we obtain the first two cases of the list.
\vskip0.2cm
\noindent

{\bf Step 2:} $n_1=1$, $n_2 \geq 2$.

Since the case $e_1=e_2=1$ will be studied in the next step in a more
general context, we start with the case $e_1 \geq 1$, $e_2=1$.
One has
\[
\varphi(2e_2,2,e_1,1,1,n_2)=-\frac{1}{2}n^2_2-\frac{1}{2}n_2<0
\qquad \forall n_2, e_1.
\]
Moreover one has $\varphi(d_1,d_2,e_1,e_2,1,n_2)>0$ when $d_i > 2e_i$,
$i=1,2$ or when $e_2\geq 2$. Finally, the case $e_1,e_2 \geq 2$ is
already studied in Proposition \ref{prodprop}, case $(a)$.
Thus the only possibility for $n_1=1$ and
$n_2 \geq 2$ is $e_1\geq 1$ and $e_2=1$ and the number of points $h$
is given by

\[
 (e_1+1)(n_2+1)-1 \geq h \geq  \frac{1}{n_2+2}\left[ \frac{(2e_1+1)(n_2+1)(n_2+2)}{2}-1\right]
\]
\noindent
 and we obtain the third case of the list.
\vskip0.2cm
{\bf Step 3:} $e_1=e_2=1$.\\
\indent In this case we write
$$\varphi(2,2,1,1,n_1,n_2)= \frac{(n_1^2n_2^2-n_1n_2^2-n_1^2n_2-2n_1^2-2n_2^2-3n_1n_2+2n_1+2n_2)}{4}
$$
Since we can suppose $n_2\geq n_1\geq 2 $ we have that $\varphi(2,2,1,1,n_1,n_2)$ is negative for 
$$
\begin{array}{cc}
n_1 & n_2\\
2 & any \\
3 & 3,4 \\
\end{array}
$$
It is easy to see that, if $d_1$ or $d_2$ are strictly greater than 2,
we do not have special effect varieties; as a matter of fact
\begin{equation*}
\begin{split}
\varphi(d_1,d_2,1,1,n_1,n_2)>& \varphi(3,2,1,1,n_1,n_2)= \\
=& \frac{(n_1+1)(n_2+1)}{12}\left[ n_1^2n_2+2n_1^2+5n_1n_2-2n_1-6n_2 \right]+n_2\\
\end{split}
\end{equation*}
and the last term is positive for the previous values of $n_1$ and $n_2$.
Thus we can conclude that $Y$ is a $2-$special effect variety in the following cases
$$
\begin{array}{lccc}
\PP^{2}\times \PP^{n_2}& (2,2) & (1,1) & m_2(n_2)<h\leq M_2(n_2) \\
\PP^{3}\times \PP^{3} &  (2,2) & (1,1) & 15 \\
\PP^{3}\times \PP^{4}& (2,2) & (1,1) & 19 \\
\end{array}$$
where $m_1,m_2,M_1$ and $M_2$ are defined by (\ref{boundcross}), i.e.
\begin{equation*}
\begin{split}
m_2(n_2):=\frac{3n_2^2+9n_2+5}{n_2+3} & \qquad M_2(n_2):= 3n_2+2.\\
\end{split}
\end{equation*}
Finally, the non-existence of $2-$special effect varieties of bidegree $(e_1,e_2)$,
with $e_1+e_2\geq 3$ in $\PP^{n_1}\times \PP^{n_2}$, $n_1,n_2\geq 2$
is a consequence of Proposition \ref{prodprop}--$(a)$ and the following claim.
\vskip0.1cm
\noindent
{\bf Claim:} $\varphi\geq 0$, for $e_1=1$, $e_2\geq2$ (resp. for $e_1\geq2$, $e_2=1$) for $n_1,n_2\geq 2$ and $d_i\geq 2e_i$, $i=1,2$;
\vskip0.1cm
In fact, one has
\begin{equation*}
\begin{split}
\varphi(2,2e_2,1,e_2,n_1,n_2)&= \frac{(n_1+1)(n_1+2)}{2}\frac{(n_2)_{(2e_2)}}{(2e_2)!}+ n_1+n_2-\\
&-\frac{(n_1+1)(n_2)_{(2e_2)}(n_1+n_2+1)}{e!}= \\
& =\frac{(n_1+1)(n_2)_{(e_2)}}{e_2!} C(e_2)+n_1+n_2
\end{split}
\end{equation*}
where
$$C(e_2)=\frac{(n_1+2)(n_2+e)_{(e_2)}}{2(e_2)_{(e_2)}}
-(n_1+n_2+1).$$
Then the proof uses the same argument of the claim in Lemma \ref{LN1},
verifying, at the end, that $C(2)\geq 0$ for $n_1,n_2\geq 2$. 
\vskip0.2cm
\noindent 
\vskip0.2cm
We analyze now the case $e_1\cdot e_2=0$.
By symmetry, it is enough to treat the case $e_2=0$
(we recall that $d_1\cdot d_2\not=0$). In this situation we have
\begin{equation}
\begin{split}\varphi(d_1,d_2,e_1,0,n_1,n_2) := & \binom{d_1+n_1}{n_1}\binom{d_2+n_2}{n_2} -\\
&-\binom{d_1-2e_1+n_1}{n_1}\binom{d_2+n_2}{n_2} -\\
&- \left[\binom{e_1+n_1}{n_1}-1\right](n_1+n_2+1).
\end{split}
\end{equation}
Since the previous function is non-decreasing in $d_1$ and $d_2$
we can start from the minimal degree $(d_1,d_2)=(2e_1,1)$ and we obtain
\begin{equation*}
\begin{split}
\varphi(2e_1,1,e_1,0,n_1,n_2) := & \binom{2e_1+n_1}{n_1}(n_2+1)- \left[\binom{e_1+n_1}{n_1}\right](n_1+n_2+1)+n_1.
\end{split}
\end{equation*}
This function is clearly increasing in $n_2$.
For the behaviour of $\varphi(2e_1,1,e_1,0,n_1,n_2)$ by respect to $n_1$ we can write
$$\varphi(2e_1,1,e_1,0,n_1,n_2) =  A(n_1,n_2,e_1)\cdot B(n_1,n_2,e_1)+n_1$$
where
$$A(n_1,n_2,e_1)=\frac{(n_1)_{(e_1)}}{(2e_1)!}$$
and
$$B(n_1,n_2,e_1)=(n_2+1)(n_1+e_1)_{(e_1)}-(n_1+n_2+1)(e_1)_{(e_1)}.$$
Both $A(n_1,n_2,e_1)$ and  $B(n_1,n_2,e_1)$ are increasing in $n_1$.
Moreover $A(n_1,n_2,e_1)\geq0$ for $n_1,n_2,e_1\geq 1$ and, by a
simple computation, one has
\begin{equation*}
\begin{split}
B(1,n_2,e_1)&=(n_2+1)(n_1+e_1)_{(e_1)}-(n_1+n_2+1)(e_1)_{(e_1)}=\\
&=(e_1+1)_{(e_1-1)}(n_2e_1-1)\geq 0 \mbox{ for } n_2,e_1\geq 1 
\end{split}
\end{equation*}
Hence $\varphi(2e_1,1,e_1,0,n_1,n_2)$ is non-decreasing in $n_1$ too
and we can study it starting from $n_1=n_2=1$. One has
$$\varphi(2e_1,1,e_1,0,1,1) = e_1 > 0 \quad \forall e_1 \geq 1.$$
Thus $Y$ is not a $2-$special effect variety for $\LL_{(d_1,d_2)}(2^h)$
if $Y$ has bidegree $(e_1,e_2)$, with $e_1\cdot e_2=0$ and $d_1\cdot d_2\not=0$.
\end{proof}

Let $Q$ be the quadric in $\PP^3$ and consider $L_1$ and $L_2$
the generators of $\mbox{Pic}(Q)$.
Denote by $\LL(a,b)$ the linear system $|aL_1+bL_2|$.

\begin{Cor}\label{suquadrica}
A curve of type $(n,1)$ (resp. of type $(1,n)$) on a quadric $Q\subset\PP^3$
is a $2-$special effect variety on $Q$ for $\LL(2n,2)(2^{2n+1})$
(resp. for $\LL(2,2n)(2^{2n+1})$).
\end{Cor}

\begin{proof}
It follows directly from the first two cases of Proposition \ref{listpapb}.
\end{proof}

We pass now to analyze the case in which $t\geq3$;
we restrict our studying to the case $e_i\not= 0$,
$\forall i=1,\dots, t$ and $n_t\geq n_{t-1}\geq \dots \geq n_1$. 

\begin{Prop}\label{supapbpc} Let $t\geq 3$.
Let $\LL:=\LL_{(d_1, \dots d_t)}(2^h)$ be a linear system of multidegree $(d_1, \dots , d_t)$,
with $d_i\not=0$ for $i=1,\dots, t$,
on $X=\PP^{n_1}\times \dots \times \PP^{n_t}$
passing through $h$ double points in general position
and let $Y$ be a divisor of multidegree $(e_1, \dots , e_t)$ on $X$
with $e_i\not=0$ for $i=1,\dots, t$.
Moreover we require that $Y$ passes simply
through the $h$ points in $\LL$.
Then $Y$ is a $2-$special effect variety on $X$ for $\LL$ only if $t=3$
and for the following values:
$$\begin{array}{cccc}
\hline 
\PP^{n_1}\times\PP^{n_2} \times \PP^{n_3} & (d_1,d_2,d_3) & (e_1,e_2,e_3) & h \\\hline
\PP^{1}\times\PP^{1} \times \PP^{1} & (2,2,2) & (1,1,1) & 7  \\
\PP^{1}\times\PP^{1} \times \PP^{2} & (2,2,2) & (1,1,1) &  11 \\
\PP^{1}\times\PP^{1} \times \PP^{3} & (2,2,2) & (1,1,1) &  15 \\
\hline
\end{array}$$ 
\end{Prop}

\begin{proof}
The result follows immediately from cases $(b)$ and $(c)$ of Proposition \ref{prodprop}. 
\end{proof}

\subsection{$\alpha-$Special effect varieties and Segre--Veronese varieties}\label{SegVer}
It is known from the literature that the speciality of a linear system
with imposed double points can be phrased in term of defectivity of certain varieties.
The reader can find more topics on these subjects,
for example, in \cite{Ciliberto}, \cite{Miranda} and \cite{Zak}.

Catalisano, Geramita and Gimigliano, in 
\cite{CaGeGi2}, \cite{CaGeGi3} and  \cite{CaGeGi4}, 
study the secant varieties of Segre--Veronese varieties, i.e. the
image of $\PP^{n_1}\times \cdots \times \PP^{n_t}$ under 
the composition of the Veronese embeddings
$\nu_{a_1} \times \cdots \times \nu_{a_t}$ followed by the Segre embedding $\rho_s$:
\[
\begin{CD}
\PP^{n_1}\times \cdots \times \PP^{n_t} @>{\nu_{a_1} \times \cdots \times \nu_{a_t}}>> \PP^{\binom{a_1+n_1}{n_1} -1}
\times \cdots \times \PP^{\binom{a_t+n_t}{n_t}-1} @>{\rho_s}>> \PP^N
\end{CD}
\]

Their results on defective Segre--Veronese varieties , or equivalenty 
on special linear systems on  $\PP^{n_1}\times \cdots \times \PP^{n_t}$
can be compared with our results on $\alpha-$special effect varieties.

A first result we mention is the following

\begin{Thm}[Theorem 2.1 in \cite{CaGeGi4}]\label{PrimoCGG}
Let $\LL:=\LL_{a_1,a_2}(2^h)$ be the linear system in
$\PP^1 \times \PP^1$ of divisors of bidegree $(a_1,a_2)$ with $h$ imposed double points.
Then $\LL$ is non-special unless
$$a_1=2d, a_2=2, d\geq 1, \mbox{ and } h=2d+1.$$
\end{Thm}

Using the first two cases of Proposition \ref{listpapb}
or Corollary \ref{suquadrica} we obtain immediately the following

\begin{Thm}\label{compare1}
The Numerical Conjecture holds
for each of the special systems listed in Theorem \ref{PrimoCGG}.
\end{Thm}

The second result we mention is related to the study of $\PP^1 \times \PP^1 \times \PP^1$. 

\begin{Thm}[Theorem 2.5 in \cite{CaGeGi4}]\label{SecondoCGG}
Let $a_1 \geq a_2 \geq a_3 \geq 1$, $\rho \in {\mathbb{N}}$.
Let $\LL:=\LL_{a_1,a_2,a_3}(2^h)$ be the linear system in
$\PP^1 \times \PP^1 \times \PP^1$ of divisors of multidegree $(a_1,a_2,a_3)$ with $h$ imposed double points.
Then $\LL$ is non-special unless
$$(a_1,a_2,a_3)=(2,2,2) \mbox{ and } h=7;$$
$$(a_1,a_2,a_3)=(2\alpha,1,1) \mbox{ and } h=2\alpha+1.$$
\end{Thm}

Once again we can try to check if there are special effect varieties
for the special systems corresponding to the defective varieties listed before.
It is easy to observe that, by numerical reasons, the second case cannot be treated
with a $2-$special effect variety.
However, using special effect configurations we can state a result as
Theorem \ref{compare1}.

\begin{Thm}\label{compare2}
The Numerical Conjecture holds
for each of the special systems listed in Theorem \ref{SecondoCGG}.
\end{Thm}

\begin{proof}
For the case $(a_1,a_2,a_3)=(2,2,2)$, $h=7$
there is a $2-$special effect variety as showed in Proposition \ref{supapbpc},
for $t=3$, $n_1=n_2=n_3=1$ and $d_i=a_i$, $i=1,2,3$. 
Let $\LL$ be the special linear system $\LL_{(2\alpha,1,1)}(2^{2\alpha+1})$
Let $Y_1$ [resp. $Y_2$] be a divisor corresponding to the system
$\LL_{(\alpha,0,1)}(1^{2\alpha+1})$ [resp. $\LL_{(\alpha,1,0)}(1^{2\alpha+1})$].
We easily compute
\vskip0.2cm
$
\begin{array}{l}
\nu(\LL)=-1 \\
\nu(Y_1)=\nu(Y_2)=0\\
\nu(\LL-Y_1)=\nu(\LL_{(\alpha,1,0)}(1^{2\alpha+1}))=0\\
\nu(\LL-Y_2)=\nu(\LL_{(\alpha,0,1)}(1^{2\alpha+1}))=0
\end{array}
$
\vskip0.2cm
\noindent
Then $Y+Y'$ is a $(1,1)-$special effect configuration for $\LL$.
\end{proof}


\end{document}